\documentclass[11pt]{amsart}
\usepackage{amssymb}
\usepackage{verbatim}

\begin{document}

\newtheorem{theorem}{Theorem}    
\newtheorem{proposition}[theorem]{Proposition}
\newtheorem{corollary}[theorem]{Corollary}
\newtheorem{lemma}[theorem]{Lemma}
\newtheorem{observation}[theorem]{Observation}
\theoremstyle{definition}
\newtheorem{definition}{Definition}
\newtheorem{remark}{Remark}
\newtheorem{conjecture}[remark]{Conjecture}
\newtheorem{question}[remark]{Question}
\newtheorem{example}{Example}
\def\theexample{\unskip}
\newtheorem{problem}{Problem}

\numberwithin{theorem}{section}
\numberwithin{definition}{section}
\numberwithin{remark}{section}
\numberwithin{equation}{section}

\def\reals{{\mathbb R}}
\def\integers{{\mathbb Z}}
\def\complex{{\mathbb C}\/}
\def\naturals{{\mathbb N}\/}
\def\distance{\operatorname{distance}\,}
\def\degree{\operatorname{degree}\,}
\def\dim{\operatorname{dimension}\,}
\def\Span{\operatorname{span}\,}
\def\ZZ{ {\mathbb Z} }
\def\eps{\varepsilon}
\def\p{\partial}
\def\rp{{ ^{-1} }}
\def\Re{\operatorname{Re\,} }
\def\Im{\operatorname{Im\,} }
\def\lt{L^2}
\def\Q{Q}

\def\scriptx{{\mathcal X}}
\def\scriptl{{\mathfrak L}}
\def\scriptd{{\mathcal D}}
\def\barl{{\bar L}}

\author{Michael Christ}
\address{
        Michael Christ\\
        Department of Mathematics\\
        University of California \\
        Berkeley, CA 94720-3840, USA}
\email{mchrist@math.berkeley.edu}
\thanks{Supported in part by NSF grants DMS-9970660 and DMS-0402160}

\date{September 10, 2004. Revised March 23, 2005.} 

\title[Sums of squares of complex vector fields]
{A remark on sums of squares \\ of complex
vector fields}


\maketitle


\section{Introduction }
Let $\{Z_j\}$ be a finite collection of 
vector fields, with smooth complex-valued coefficients,
defined in an open subset $U$ of Euclidean space.
Let $Z_j^*$ be the formal adjoint of $Z_j$,
with respect to the Hilbert space structure $L^2$ associated
to some measure with a smooth nonvanishing density.
Consider the operator 
$\scriptl = \sum_j Z_j^*Z_j$, which we shall refer to as a sum of squares.
$\scriptl$ is said to be hypoelliptic in $U$ if 
for any open subset $V\subset U$
and any distribution $u\in \scriptd'(V)$ such that
$\scriptl(u)\in C^\infty(V)$, necessarily $u\in C^\infty(V)$.

Assume throughout this paragraph only that all vector fields are real.
Then a well-known sufficient condition for hypoellipticity
is the bracket condition of H\"ormander, that the Lie algebra
generated by $\{Z_j\}$ should span the tangent space to
$U$ at each of its points. This condition 
ensures, and is equivalent
to, the condition that $\scriptl$ is subelliptic in the sense
that for any relatively compact open subset $V\Subset U$,
there exist $\eps>0$ and $C<\infty$ such that
for all $u\in C_0^2(V)$,
\begin{equation} \label{subelliptic}
\|u\|_{H^\eps} \le C\|\scriptl u\|_{H^0} + C\|u\|_{H^0}.
\end{equation}
This can be equivalently reformulated as
\begin{equation} \label{quadformsubelliptic}
\|u\|_{H^\eps}^2 \le CQ(u,u) + C\|u\|_{H^0}
\end{equation}
where $Q(u,u) = \sum_j \|Z_j u\|_{H^0}^2$.
Subellipticity in turn implies hypoellipticity for sums of squares operators.
However, $\scriptl$ is sometimes hypoelliptic
without satisfying the bracket condition. See for instance \cite{christdegenregime}
and the references cited there.

Henceforth we allow vector fields to be complex.
Weaker inequalities than \eqref{subelliptic} are then conceivable.
\begin{definition} \label{lossdefn}
We say that $\scriptl$ loses at most finitely many derivatives
in any open set $U$
if for every $V\Subset U$
there exist $s>-\infty$, $t<+\infty$ and $s'<s$ such that for all $u\in C^\infty_0(V)$,
\begin{equation} \label{lossy}
\|u\|_{H^s}\le C\|\scriptl u\|_{H^t}
+C\|u\|_{H^{s'}}.
\end{equation}
We say that it loses derivatives\footnote{Many authors use these words differently,
considering any nonelliptic operator to lose derivatives; for an operator of order
$m$, losing derivatives in our language corresponds essentially to losing
more than $m$ derivatives in that alternative language.}
if for any $t$, no such inequality holds with $s=t$.
\end{definition}
\noindent
This usage is not universally accepted, and will be discussed further in
\S\ref{section:whatlossmeans} below.

For complex vector fields
the bracket condition still makes sense, and subellipticity in the sense
\eqref{quadformsubelliptic} continues to imply hypoellipticity.
Siu has asked whether the bracket condition continues to imply subellipticity in this sense
for complex fields.
Kohn \cite{kohn} has answered this in the negative\footnote{But has shown that it
does imply subellipticity if $\{Z_j\}$ together with their
brackets with only two factors suffice to span the tangent space.}, 
and has gone further by establishing examples which simultaneously
(i) satisfy the bracket hypothesis, (ii) not only fail to be subelliptic but actually
lose derivatives, yet (iii) are nonetheless hypoelliptic.
This note is a comment on \cite{kohn}, showing that even the weaker property of
hypoellipticity can fail, for complex vector fields satisfying the bracket condition.

Earlier, Heller \cite{heller} had studied the hypoellipticity
(and analytic hypoellipticity) of left-invariant differential
operators of arbitrary order on the Heisenberg group, subject
to a hypothesis of transversal ellipticity. 
He showed that such an operator is ($C^\infty$ and $C^\omega$) hypoelliptic 
whenever it loses at most finitely many derivatives, 
and he gave an example of a fourth order operator\footnote{ 
Namely $\square_b^2+X$ on ${\mathbb H}^1$.}
which does lose derivatives, yet is hypoelliptic.
This extended an analysis of 
Stein \cite{stein}, who had proved hypoellipticity (as well as analytic hypoellipticity) 
for certain second order operators\footnote{ Such as $\square_b+1$.}
which do not gain derivatives, but do not actually lose them either. 

Hypoellipticity with loss of derivatives is a delicate
matter, because estimates without any gain
in regularity are inevitably quite unstable. Any analysis
of hypoellipticity must involve deformation of $\scriptl$, for instance
via the introduction of some type of cutoff operators, 
potentially destroying the estimates \eqref{lossy}.

From this point of view our main result is not surprising:
\begin{proposition} \label{prop:imprecise}
There exist finite families of complex vector fields
$Z_j$ with $C^\infty$ coefficients
which satisfy the bracket condition and lose at most finitely
many derivatives in the sense \eqref{lossy},
but for which $\sum_j Z_j^*Z_j$ fails to be $C^\infty$ hypoelliptic.
\end{proposition}

To describe these consider $\reals^3$ with coordinates
$(x,t,s)$. We consider always the Hilbert space $L^2(\reals^3)$ associated to
Lebesgue measure in these coordinates.
Define
\begin{equation}
\barl = \p_x-ix\p_t,\qquad L = \p_x+ix\p_t.
\end{equation}
Fix an integer $k\ge 1$ and define
$Z_1 = \barl$, $Z_2 = x^k L$, and $Z_3 = \p_s$.
Here $\p_x = \frac{d}{dx}$, with no factor of $\sqrt{-1}$, and so forth.
Proposition~\ref{prop:imprecise} can now be more precisely restated.

\begin{proposition}
Let $k$ be any positive integer.
The complex vector fields $Z_1,Z_2,Z_3$ satisfy the bracket
condition at each point of $\reals^3$, and in any bounded open
set $V\subset\reals^3$, the operator 
$\scriptl = \sum_{j=1}^3 Z_j^*Z_j$ 
loses at most finitely many derivatives.
Nonetheless, $\scriptl$ is not $C^\infty$ hypoelliptic in 
any neighborhood of the origin.
\end{proposition}
When $k=1$, $\scriptl$ actually satisfies \eqref{lossy} with
$s=0$, that is, it does not lose derivatives; yet it fails to
be hypoelliptic.

$Z_1,Z_2$ can also be regarded as vector fields in $\reals^2$ 
rather than in $\reals^3$.
The operator $\scriptl_{\reals^2} = Z_1^*Z_1+Z_2^*Z_2$ in $\reals^2$
is then a simplified version of Kohn's examples, and 
can be shown to be hypoelliptic although we will not do so here.
Adding the extra variable $s$ and the extra term $-\p_s^2$ to create $\scriptl$
destroys hypoellipticity, due to propagation
of singularities along curves such as $\{(0,0,s)\}$.

Our example is closely analogous
to two well-known examples concerning $C^\infty$ and analytic
hypoellipticity \cite{baouendigoul}, \cite{kusuokastroock}. 
Firstly, the operator
$-\p_x^2-x^2\p_t^2$
is analytic hypoelliptic in $\reals^2$, whereas
$-\p_x^2-x^2\p_t^2-\p_s^2$
fails to be analytic hypoelliptic in $\reals^3$.
Secondly,
consider a $C^\infty$ function $a:\reals^1\to\reals$ 
such that $a(x)=0$ if and only if $x=0$.
Then
$-\p_x^2-a(x)^2\p_t^2$
is always $C^\infty$ hypoelliptic in $\reals^2$,
while
$-\p_x^2-a(x)^2\p_t^2 -\p_s^2$
may or may not be hypoelliptic in $\reals^3$,
depending on the rate at which $a(x)$ tends to zero
as $x\to 0$.
For an attempt to place these examples in perspective
see \cite{christdegenregime}, \cite{christspec}.

The author is indebted to Joe Kohn for stimulating discussions.

\section{Spectral analysis of certain ODEs}
For $\tau\in\reals^+$
consider the ordinary differential operators
\begin{equation}
P_\tau = 
-(\p_x-x\tau)(\p_x+x\tau)
- (\p_x+x\tau)x^{2k}(\p_x-x\tau).
\end{equation}
These are obtained by separation of variables;
\[\scriptl_{\reals^2}(e^{i\tau t}f(x)) = e^{i\tau t}P_\tau f(x).\]
$P_\tau$ is formally selfadjoint on $\lt(\reals)$ with
respect to Lebesgue measure, and is nonnegative. 

For any $\tau>0$, $P_\tau $ is unitarily equivalent,
via the change of variables $y = \tau^{1/2}x$
and substitution $F(y) = \tau^{-1/4}f(x)$, to
$\tau\Q_\tau$ where
\begin{equation}
\Q_\tau = 
-(\p_y-y)(\p_y+y)
- \tau^{-k}(\p_y+y)y^{2k}(\p_y-y).
\end{equation}
Setting $g(y) = e^{-y^2/2}$ we have
\begin{equation} \label{odeapriori}
\langle \Q_\tau g,\,g\rangle
 = \tau^{-k}\|y^k(\p_y-y) e^{-y^2/2}\|_{\lt}^2
= c\tau^{-k}.
\end{equation}

Conversely we claim that for all $\tau\ge 1$ and all $f\in C^2_0(\reals)$,
\begin{equation} \label{globalelliptic}
\tau^{-k}\|f\|_{L^2}^2
+ 
\tau^{-k}\|yf\|_{L^2}^2
+ 
\tau^{-k}\|\p_y f\|_{L^2}^2
\le C\langle \Q_\tau f,\,f\rangle.
\end{equation}
Indeed,
\begin{align*}
\langle \Q_\tau f,\,f\rangle
&=
\tau^{-k}\|y^{k} (\p_y-y)f\|_{\lt}^2
+ \|(\p_y+y)f\|_{\lt}^2
\\
&\ge \tau^{-k}\int_{|y|\ge 1} |(\p_y-y)f|^2\,dy
+ \int |(\p_y+y)f|^2\,dy
\\
&\ge \tau^{-k}\int_{|y|\ge 1} |yf(y)|^2\,dy
+ \int_{|y|\le 1} |((\p_y+y)f|^2\,dy
\\
&\ge c\tau^{-k}\|f\|_{\lt}^2,
\end{align*}
and \eqref{globalelliptic} follows from this together with
the majorization
$\langle \Q_\tau f,\,f\rangle \ge \|(\p_y+y)f\|^2_{\lt}$.

It follows readily that the $\lt$ closure of $\Q_\tau$ is
selfadjoint and has discrete spectrum, and that every
eigenfunction of $\Q_\tau$ belongs to the Schwartz space. 
Define $\lambda(\tau)$
to be the lowest eigenvalue of
$\Q_\tau$.
By \eqref{odeapriori} and \eqref{globalelliptic}, there exist $0<c<c'<\infty$ such that 
\begin{equation} \label{eigenvaluebound}
c'\tau^{-k}\le \lambda(\tau) \le c\tau^{-k}
\qquad\forall\,\,\tau\in[1,\infty).
\end{equation}

Let $\psi_\tau\in\lt(\reals)$ be an eigenfunction of $\Q_\tau$
with eigenvalue $\lambda(\tau)$, normalized so
that $\|\psi_\tau\|_{\lt(\reals)}=1$.
We claim that 
 $\|y\psi_\tau\|_{\lt}$ and $\|\p_y\psi_\tau\|_{\lt}$
are bounded above, uniformly in $\tau$ for all $\tau\ge 1$.
To prove this, decompose $\psi_\tau = ah_0+g$
where $h_0(y) = e^{-y^2/2}$, $a\in\complex$,
and $g\perp h_0$.
Since $\p_y+y$ annihilates $h_0$, since $\|(\p_y+y)\psi_\tau\|^2_{\lt}\le
\langle \Q_\tau \psi_\tau,\,\psi_\tau\rangle$,
and since $\|g\|_{\lt}\lesssim\|(\p_y+y)g\|_{\lt}$,
it follows that for large $\tau$ one has
$\|g\|_{\lt}\lesssim\lambda(\tau)\ll 1$, 
and consequently
$|a|\sim 1$.
Since $\|yg\|_{\lt}+\|\p_y g\|_{\lt}\lesssim \|(\p_y+y)g\|_{\lt} + \|g\|_{\lt}$
for all functions $g$ orthogonal to $h_0$,
and since $h_0$ is a Schwartz function and is independent of $\tau$, the 
claim follows.

From this we conclude firstly that $\|\psi_\tau\|_{L^\infty(\reals)}$
is bounded above, uniformly for all $\tau\ge 1$.
Secondly there exists $B<\infty$ such that 
\begin{equation} \label{psibignear0}
\sup_{|y|\le B} |\psi_\tau(y)|\ge B\rp
\end{equation}
uniformly for all $\tau\ge 1$.

\section{Conclusion of proof}
Consider the family of functions $u_\tau$ defined for $\tau\in[1,\infty)$ by
\begin{equation}
u_\tau(x,t,s)
= e^{i\tau t} e^{\sigma(\tau) s} 
\psi_\tau(\tau^{1/2}x)
\end{equation}
where $\sigma(\tau)>0$ is the positive solution of 
$\sigma^2 = \tau\lambda(\tau)$.
Then $\scriptl u_\tau\equiv 0$ in $\reals^3$.
By \eqref{eigenvaluebound}, $\sigma(\tau) = O(\tau^{(1-k)/2})$; in particular,
$\sigma(\tau)$ is uniformly bounded as $\tau\to+\infty$.

As is well known, hypoellipticity implies certain inequalities
via the Baire category theorem.
If $\scriptl$ were hypoelliptic,
then for any open sets $V\Subset V'$ and any $N\in\naturals$
there would exist $C,M<\infty$ 
such that for all $u\in C^\infty(V')$,
\begin{equation} \label{baire}
\|u\|_{C^N(V)}
\le C\|\scriptl u\|_{C^M(V')}
+ C\|u\|_{C^0(V')}.
\end{equation}

Fix $V\Subset V'\Subset\reals^3$ with $0\in V$.
Consider the inequality \eqref{baire}
for $u_\tau$, for large positive $\tau$.
By \eqref{psibignear0}, for all sufficiently large $\tau$
we have
\begin{equation}
\|\p_t u_\tau\|_{C^0(V)}
\ge c\tau.
\end{equation}
On the other hand $\scriptl u_\tau\equiv 0$,
while the uniform boundedness of $\psi_\tau$ in $L^\infty$ implies that
\begin{equation}
\|u_\tau\|_{C^0(V')}
\le C\|\psi_\tau\|_{C^0(\reals)} e^{C\sigma(\tau)}
\le C'e^{C\sigma(\tau)};
\end{equation}
the factor $e^{\sigma(\tau) s}$ is $O(e^{C\sigma(\tau)})$ because $V'$ is a bounded set.
Since $\sigma(\tau) = O(\tau^{(1-k)/2})$ remains bounded
as $\tau\to\infty$, $\|u_\tau\|_{C^0(V')}$ likewise
remains uniformly bounded. 
Thus \eqref{baire} fails to hold for $N=1$.
\qed

\section{On loss of derivatives} \label{section:whatlossmeans}

Definition~\ref{lossdefn} is only one possible notion of loss of derivatives.
A more common notion, as Kohn has pointed out,
is essentially this: 
$\scriptl$ is said to lose  at least $\delta$ derivatives
in an open set $U$ if there exist an open subset $V\subset U$, an exponent $s$, 
and a distribution $u\in\scriptd'(V)$
such that $\scriptl u\in H^s_{\rm loc}(V)$,
yet $u\notin H^t_{\rm loc}(V)$ for any $t>s-\delta$.
$\scriptl$ is then said to lose derivatives 
if it loses at least $\delta$ derivatives for some $\delta>0$.
It is thus formally conceivable that an operator could lose at most 
a certain number of derivatives in the sense of Definition~\ref{lossdefn},
yet lose more derivatives, or even infinitely many, in this alternative sense.

A global inequality of the form \eqref{lossy} expresses a very weak property of an
operator.
Hypoellipticity amounts to having a family of inequalities that are stronger
in two ways, incorporating both (i) spatial localization and (ii) a type
of localization (expressed by weighted $L^2$ inequalities) with respect
to frequency variables in phase space.
An inequality corresponding to an implication
$\scriptl u\in H^s_{\rm loc}\Rightarrow u\in H^t_{\rm loc}$
expresses one of these two types of localization, but not the other.
We regard such an inequality as expressing a type of partial hypoellipticity,
whereas \eqref{lossy} is a minimal {\it a priori} inequality
involving no localization. 
We note that such inequalities, with $\scriptl$
replaced by its transpose, are fundamental to the theory of local solvability.
\eqref{lossy} appears at one extreme of a (partially ordered) spectrum 
of possible inequalities,
with hypoellipticity lying at the opposite end of the spectrum
and the notion of loss discussed in the preceding paragraph lying
somewhere in between.

Other variants formulated in terms of the quadratic form
$Q(u,u)=\sum_j\|Z_j u\|_{H^0}^2$, rather than some norm of $\scriptl u$,
are also reasonable.

\end{document}